\documentclass[11pt]{amsart}\usepackage{latexsym}\usepackage{amsfonts}\usepackage{amsmath}\usepackage{amssymb}\usepackage{amscd}\usepackage[cm]{fullpage}
\usepackage{color,CJK}
\newtheorem{thm}{Theorem}[section]
\newtheorem{prop}{Proposition}[section]

\newtheorem{remark}{Remark}[section]

\begin{document}
\title[Curvature decomposition]{Curvature decompositions on Einstein four-manifolds}
\author{Peng Wu}
\address{Shanghai Center for Mathematical Sciences, Fudan University, Shanghai 200433, China} \email{wupenguin@fudan.edu.cn}\thanks{}
\dedicatory{}\date{\today}
\begin{abstract}
For Einstein four-manifolds with positive scalar curvature, we derive relations among various positivity conditions on the curvature tensor, some of which are of great importance in the study of the Ricci flow. These relations suggest possible new ideas to study the well-known rigidity conjecture for positively curved Einstein four-manifolds.
\end{abstract}
\maketitle

\section{Introduction}

A Riemannian metric is called an Einstein metric if $\mathrm{Ric}=\lambda g$ for some $\lambda\in\mathbb{R}$. A central problem in differential geometry is to study the existence, rigidity, and moduli space of Einstein metrics. In dimension four, a well-known conjecture states that Einstein four-manifolds with positive sectional curvature are isometric to $(S^4,g_0)$ or $(\mathrm{C}P^2,g_{FS})$. Many authors have made important progress on this conjecture, cf. Berger \cite{berger61}, Derdzinski \cite{Der}, Hitchin \cite{Besse}, Gursky and LeBrun \cite{GL}, Yang \cite{Yang}, and Costa \cite{Costa}. Curvature decompositions are basic tools to understand the structure of the curvature tensor. The three curvature decompositions on Einstein four-manifolds: the standard curvature decomposition, the duality curvature decomposition, and the Berger curvature decomposition, are essential in these works.

The positivity of the curvature operator is of great importance in the study of the Ricci flow. Recall that a curvature operator $\mathfrak{R}$ is $k$-positive ($k$-nonnegative), if the sum of its $k$ smallest eigenvalues is positive (nonnegative). In a pioneering work, Hamilton \cite{H4} proved that the space of positive curvature operator is preserved along the Ricci flow, and compact four-manifolds with positive curvature operator are diffeomorphic to spherical space forms. Chen \cite{Chen} later relaxed Hamilton's condition to 2-positive curvature operator. In a recent breakthrough, B\"ohm and Wilking \cite{BW} proved that compact $n$-dimensional manifolds with 2-positive curvature operator are diffeomorphic to spherical space forms. Unfortunately, as B\"ohm and Wilking \cite{BW} pointed out, the space of 3-positive curvature operator is not preserved along the Ricci flow. However since the curvature operator of $(\mathrm{C}P^2,g_{FS})$ is 3-positive, it is natural to study the rigidity of Einstein four-manifolds with 3-positive or 4-positive curvature operator. As the first step, %in the author's doctoral thesis \cite{Wuthesis},
we investigate the relationship among $k$-positive curvature operator, positive sectional curvature, and positive isotropic curvature (see Section 2 for the definition).
\begin{thm} \label{thm}
Let $(M^4, g)$ be an Einstein four-manifold with $\mathrm{Ric}=\lambda g$, $\lambda>0$.
\begin{enumerate}
\item $\mathfrak{R}$ is 2-positive if and only if the isotropic curvature is positive.
\item If $K>\frac{\lambda}{12}$, then $\mathfrak{R}$ is 3-positive; if $\mathfrak{R}$ is 3-positive, then $K>\frac{\lambda}{30}$.
\item $\mathfrak{R}$ is 4-positive if and only if $K<\lambda$, which implies $K>(4-\sqrt{17})\lambda$.
\end{enumerate}
\end{thm}

The rigidity of Einstein manifolds with positive curvature operator and positive isotropic curvature have been studied by Tachibana \cite{Tachi} and Brendle \cite{Brendle}. Tachibana \cite{Tachi} proved that Einstein manifolds with positive curvature operator are isometric to spherical space forms. Brendle \cite{Brendle} proved that Einstein manifolds with positive isotropic curvature are isometric to spherical space forms.

\begin{remark}
If furthermore the metric is Hermitian, then it follows from a direct calculation that $4$-positive curvature operator is equivalent to positive orthogonal bisectional curvature.
\end{remark}

The basic idea of the proof, motivated by the work of Brendle \cite{Brendle}, is to apply the maximum principle to an equation of the curvature tensor, and reduce the problem to constrained optimizations. The new ingredient in the proof is to combine an analog of Brendle's argument \cite{Brendle} and the Berger curvature decomposition.

Notice that $K>\frac{\lambda}{12}$ implies $K<\frac{5\lambda}{6}$. Using the same argument as in Theorem \ref{thm}, we can show that a slightly smaller upper bound also implies 3-positive curvature operator,
\begin{prop} \label{remark1}
Let $(M, g)$ be an Einstein four-manifold with $\mathrm{Ric}=\lambda g$, $\lambda>0$. If $K<\frac{14-\sqrt{19}}{12}\lambda\approx(\frac{5}{6}-\frac{3}{100})\lambda$, %\approx 0.8034$
then $\mathfrak{R}$ is 3-positive.
\end{prop}

The proof of Theorem \ref{thm} shows that on Einstein four-manifolds, the upper bound and lower bound of the sectional curvature are asymmetric. For simplicity, we assume $\lambda=1$. On one hand, $K\geq\delta$ implies $K\leq1-2\delta$. For example $\delta=\frac{1}{6}$ for $(\mathrm{C}P^2,g_{FS})$. On the other hand, $K\leq\delta$ (naively) implies $K\geq1-2\delta$. However by our argument, the lower bound can be made much larger than $1-2\delta$. For example, 4-nonnegative curvature operator (equivalently $K\leq1$) implies $K\geq-1$, but from Theorem \ref{thm} we can make $K\geq4-\sqrt{17}$. This suggests that $K<1$ may be equivalent to $K>0$. Half Weyl curvature and half curvature operator have a similar asymmetric property. We denote eigenvalues of $W^{\pm}$ by $\lambda_1\leq\lambda_2\leq\lambda_3$. Notice that $-2\lambda_3\leq\lambda_1\leq-\frac{1}{2}\lambda_3$ since $W^{\pm}$ is traceless.
\begin{prop} \label{prop2}
Let $(M, g)$ be an Einstein four-manifold with $\mathrm{Ric}=g$. Suppose the minimum of $\lambda_1$ is achieved at $p$. Then $\lambda_1(p)\geq\frac{1}{2}(2\lambda_3+1-\sqrt{12\lambda_3^2+4\lambda_3+1}\,)(p)>(1-\sqrt{3}\,)\lambda_3(p)$.
\end{prop}

The proof of Theorem \ref{thm} also provides an alternative proof of the Weitzenb\"ock formula for Einstein metrics on four-manifolds by Derdzinski \cite{Der}. Moreover the alternative proof directly extends from Einstein metrics on four-manifolds to ``Einstein metrics" on four-dimensional smooth metric measure spaces, including gradient Ricci solitons, quasi-Einstein metrics, etc (see \cite{Wu1,Wu2} for details).

For readers' convenience, we now provide the following table of curvature conditions for Einstein metrics on four-manifolds,
\begin{center}
\begin{tabular}{|cclcrc|}
%\begin{tabular*}{\textwidth}{l @{\extracolsep{\fill}} clcrc}
\hline
$\mathfrak{R}$\text{ positive} $\Rightarrow$ &$\mathfrak{R}$\text{ 2-positive} $\Rightarrow$ & $K>\frac{1}{12}$ $\Rightarrow$ & $\mathfrak{R}$\text{ 3-positive} $\Rightarrow$ & $K>\frac{1}{30}$ &$\Rightarrow$\ \ $K>0$ \\
&$\Updownarrow$ & & & & $\Downarrow$\\
&\text{PIC} & & & & $\mathfrak{R}$\text{ 4-positive} \\
%\hline
& $\Downarrow$ & & & & $\Updownarrow$\\
\text{half 2-positive}&$ \Leftrightarrow$\ \text{half PIC}\quad \ \ & & & & $K<1$ \\
& $\Downarrow$ & & & & $\Downarrow$ \\
& \text{conf. half PIC} & & & $R>0\ \Leftrightarrow$ &$\mathfrak{R}$\text{ 6-positive}\\ %$K>\frac{(7-\sqrt{105})R}{112}$\\
\hline
\end{tabular}\\
Table 1. Curvature table for Einstein metrics on four-manifolds (Wu \cite{Wuthesis,Wu2}).
\end{center}
\vspace{-0.1cm}
Here $R$ is the scalar curvature; PIC denotes positive isotropic curvature; half PIC means PIC for orthonormal four-frame of a fixed orientation; and conformally half PIC means that there is a metric with half PIC in the conformal class of the Einstein metric.

From above relations, it is natural to ask the following questions on Einstein four-manifolds.
\begin{enumerate}
\item If the curvature operator is 3-positive, is $(M,g)$ isometric to $(S^4,g_0)$ or $(\mathrm{C}P^2,g_{FS})$?

\item If the sectional curvature is positive, is the curvature operator 3-positive?

\item If the curvature operator is 4-positive, is the sectional curvature positive?
\end{enumerate}

Question (1) is answered in a sequel \cite{Wu1} to the author's thesis, yet the rest two remain open.

\textbf{Acknowledgement.} This paper is based on a part of the author's Ph.D. thesis at University of California, Santa Barbara in 2012. The author thanks his advisors Professors Xianzhe Dai and Guofang Wei for their guidance, encouragement, and constant support. He thanks Professors Jeffrey Case and Jingrun Chen for helpful discussions. The author thanks the anonymous referee for many helpful suggestions. The author was partially supported by China Recruitment Program of Global Experts.

\section{Proof of Results}

We first summarize the three curvature decompositions on Einstein four-manifolds: the standard curvature decomposition, the duality curvature decomposition, and the Berger curvature decomposition.

On a Riemannian manifold $(M^n, g)$, the standard curvature decomposition is
\begin{equation*}
\mathrm{Rm}=W+\frac{1}{n-2}\mathrm{Ric}\odot g-\frac{R}{2(n-1)(n-2)}g\odot g
=W+\frac{1}{n-2}\overset{\circ}{\mathrm{Ric}}\odot g+\frac{R}{2n(n-1)}g\odot g.
\end{equation*}

On an oriented four-manifold $(M^4, g)$, the Hodge star operator $\star: \wedge^2 TM \rightarrow \wedge^2 TM$ induces a natural decomposition of the vector bundle of 2-forms $\wedge^2 TM$,
\begin{equation*}
\wedge^2 TM =\wedge^+ M \oplus \wedge^- M,
\end{equation*}
where $\wedge^{\pm} M$ are eigenspaces of $\pm 1$ respectively, sections of which are called self-dual, anti-self-dual 2-forms. It further induces a decomposition for the curvature operator $\mathfrak{R}: \wedge^2 TM \rightarrow \wedge^2 TM$,
\begin{equation*}
\mathfrak{R} =\left( \begin{array}{cc}
\frac{R}{12}g+W^+ & \overset{\circ}{\mathrm{Ric}} \\
\overset{\circ}{\mathrm{Ric}} & \frac{R}{12}g+W^-
\end{array} \right) ,
\end{equation*}
where $\overset{\circ}{\mathrm{Ric}}$ is the traceless Ricci curvature, $R$ is the scalar curvature. In particular if $(M^4,g)$ is an Einstein manifold, then
\begin{equation}
\label{dual} \mathfrak{R}=\left( \begin{array}{cc}
\frac{R}{12}g+W^+ & 0 \\
0 & \frac{R}{12}g+W^-
\end{array} \right) \triangleq\left( \begin{array}{cc}
\mathfrak{R}^+ & 0 \\
0 & \mathfrak{R}^-
\end{array} \right).
\end{equation}

In \cite{berger61}, Berger discovered another curvature decomposition for Einstein four-manifolds (see also Singer and Thorpe \cite{ST}),
\begin{prop} \label{Bergerdecomposition}
\label{prop2.3} Let $(M,\ g)$ be an Einstein four-manifold with $\mathrm{Ric}=\lambda g$. For any $p\in M$, there exists an orthonormal basis $\{e_i\}_{1\leq i\leq 4}$ of $T_p M$, such that relative to the corresponding basis $\{e_i\wedge e_j\}_{1\leq i<j\leq 4}$ of $\wedge^2 T_pM$, $\mathfrak{R}$ takes the form
\begin{equation}
\label{abba} \mathfrak{R}=\left( \begin{array}{cc}
A & B\\
B & A
\end{array}\right),
\end{equation}
where $A=\mathrm{diag}\{a_1,\ a_2,\ a_3\}$, $B=\mathrm{diag}\{b_1,\ b_2,\ b_3\}$ satisfying the following properties,\\
$\mathrm{(1)}$. $a_1=K(e_1, e_2)=K(e_3, e_4)=\min\{K(\sigma): \sigma\in\wedge^2 T_p M,\|\sigma\|=1\}$,\\
\quad\quad $a_3=K(e_1,e_4)=K(e_2,e_3)=\max\{K(\sigma):\sigma\in\wedge^2 T_p M, \|\sigma\|=1\}$,\\
\quad\quad $a_2=K(e_1, e_3)=K(e_2, e_4)$, and $a_1+a_2+a_3=\lambda$;\\
$\mathrm{(2)}$. $b_1=R_{1234},\ b_2=R_{1342},\ b_3=R_{1423}$;\\
$\mathrm{(3)}$. $|b_2-b_1|\leq a_2-a_1,\ |b_3-b_1|\leq a_3-a_1,\ |b_3-b_2|\leq a_3-a_2$.
\end{prop}

Diagonalizing the matrix in the Berger curvature decomposition, we get eigenvalues of curvature operator $\mathfrak{R}$ and half curvature operators $\mathfrak{R}^{\pm}$ in the following order,
\begin{equation}\label{eigenvalue}
\begin{cases}
& a_1+b_1 \leq a_2+b_2 \leq a_3+b_3,\\
& a_1-b_1 \leq a_2-b_2 \leq a_3-b_3.
\end{cases}
\end{equation}

Therefore by the Berger curvature decomposition, we have

$\bullet$ positive sectional curvature is equivalent to $(a_1+b_1)+(a_1-b_1)>0$, that is, the sum of the smallest eigenvalues of  $\mathfrak{R}^+$ and $\mathfrak{R}^-$ is positive;

$\bullet$ 2-positive curvature operator is equivalent to $(a_1+a_2)\pm(b_1+b_2)>0$ and $a_1>0$;

$\bullet$ positive isotropic curvature implies $(a_1+a_2)\pm(b_1+b_2)>0$;

$\bullet$ 3-positive curvature operator is equivalent to $2a_1+a_2\pm b_2>0$;

$\bullet$ 4-positive curvature operator is equivalent to $a_1+a_2>0$ and $1+(a_1\pm b_1)>0$.

Recall that $(M,g)$ is said to have positive isotropic curvature \cite{MM}, if for any orthonormal four-frame $\{e_i,e_j,e_k,e_l\}$, the curvature tensor satisfies
\begin{equation*}
\begin{split}
R_{ikik}+R_{ilil}+R_{jkjk}+R_{jljl}>2R_{ijkl}.
\end{split}
\end{equation*}

In fact the Berger curvature decomposition corresponds to a special duality curvature decomposition, because eigenvectors of $a_i\pm b_i$ are self-dual and anti-self-dual $2$-forms, respectively.

\textbf{Proof} of Theorem \ref{thm}. Without loss of generality we assume $\lambda=1$. We start with some simple observations. %(4) and (8) are obvious.
It is well known that 2-positive curvature operator implies positive isotropic curvature. By the Berger curvature decomposition, we have
\begin{equation*}
\begin{split}
&a_1-a_2\leq b_2-b_1\leq a_2-a_1,\\
&a_2-a_3\leq b_2-b_3\leq a_3-a_2,
\end{split}
\end{equation*}
taking the sum we get $|b_2|\leq\frac{1}{3}(a_3-a_1)$. If $a_1>\frac{1}{12}$, then
\begin{equation*}
2a_1+a_2-|b_2|\geq 2a_1+a_2-\frac{1}{3}(a_3-a_1)\geq 4a_1-\frac{1}{3}> 0,
\end{equation*}
so $\mathfrak{R}$ is 3-positive. If $\mathfrak{R}$ is 4-positive, it is obvious that $a_1+a_2>0$, so $K<1$.

\

Recall that for Einstein manifolds (see Hamilton \cite{H3}),
\begin{equation} \label{Hamiltonequation}
\Delta R(e_i,e_j,e_k,e_l) +2(B_{ijkl}-B_{ijlk}+B_{ikjl}-B_{iljk})
=2R_{ijkl},
\end{equation}
where $B_{ijkl}=g^{mn}g^{pq}R_{imjp}R_{knlq}$. Applying the Berger curvature decomposition, we get
\begin{equation*}
\begin{cases}
\Delta R(e_1,e_2,e_1,e_2) +2(a_1^2+b_1^2+2a_2a_3+2b_2b_3) = 2a_1,&\\
\Delta R(e_1,e_3,e_1,e_3) +2(a_2^2+b_2^2+2a_1a_3+2b_1b_3) = 2a_2,&\\
\Delta R(e_1,e_4,e_1,e_4) +2(a_3^2+b_3^2+2a_1a_2+2b_1b_2) = 2a_3.&
\end{cases}
\end{equation*}

Suppose that the minimum of the sectional curvature is attained at $p$ by the tangent plane spanned by $\{e_1,e_2\}$.  Since $2\min K=\min(\mathfrak{R}^+ +\mathfrak{R}^-)$, for any $v\in T_p M$ and the geodesic $\gamma(t)$ with $\gamma(0)=p,\ \gamma'(0)=v$, let $\{e_1,\ e_2,\ e_3,\ e_4\}$ be a parallel orthornormal frame along $\gamma(t)$, then we have
\begin{equation*}
\begin{split}
(D^2_{v,v}R)(e_1,e_2,e_1,e_2)(p)=D^2_{v,v}(R(e_1,e_2,e_1,e_2))(p)\geq0.
\end{split}
\end{equation*}
Taking the trace we have $(\Delta R)(e_1,e_2,e_1,e_2)(p) \geq 0$, therefore at $p$ we get
\begin{equation}  \label{Hamilton}
\boxed{a_1^2+b_1^2+2(a_2a_3+b_2b_3)\leq a_1.}
\end{equation}
\vspace{0.1cm}

\textit{First we prove 2-positive curvature operator is equivalent to positive isotropic curvature.} It suffices to show that $(a_1+a_2)\pm(b_1+b_2)>0$ implies $a_1>0$. In fact if $(a_1+a_2)\pm (b_1+b_2)>0$, then
\begin{equation*}
a_2\pm b_2>0,\quad\quad a_3\pm b_3>0.
\end{equation*}
Therefore by (\ref{Hamilton}), we have
\begin{equation*}
\begin{split}
a_1(p)\geq&a_1^2+b_1^2+2(a_2a_3+b_2b_3)>a_1^2+b_1^2\geq 0.
\end{split}
\end{equation*}

\textit{Next we prove 3-positive curvature operator implies positive sectional curvature.} If $\mathfrak{R}$ is 3-positive, then
\begin{equation*}
a_2\pm b_2>-2a_1, \quad\quad a_3\pm b_3>-2a_1.
\end{equation*}
Assuming that $a_1(p)\leq 0$, then $a_2\pm b_2>0$ and $a_3\pm b_3>0$, hence we have
\begin{equation*}
\begin{split}
a_1(p)\geq& a_1^2+b_1^2+2(a_2a_3+b_2b_3)>a_1^2+b_1^2\geq 0,
\end{split}
\end{equation*}
which contradicts to (\ref{Hamilton}). Therefore $a_1(p)>0$, i.e., $(M,\ g)$ has positive sectional curvature.

\vspace{0.3cm}

\textit{Next we derive a lower bound for the sectional curvature when $\mathfrak{R}$ is 3-positive.} Let $a_2(p)=ka_1(p),\ k\geq 1$. If $b_2b_3\geq 0$, then from (\ref{Hamilton}),
$$a_1 \geq a_1^2+2a_2a_3\geq a_1^2+2a_1(1-2a_1)=2a_1-3a_1^2,$$
which implies that $a_1=1/3$.

If $b_2b_3<0$, without loss of generality, we assume $b_2<0, \ b_3>0$. On one hand, by 3-positivity of the curvature operator, $|b_2|<a_2+2a_1=(k+2)a_1$, so we get
\begin{equation*}
\begin{split}
b_1^2+2b_2b_3=&b_2^2+b_3^2+4b_2b_3
=(b_3+2b_2)^2-3b_2^2
>-3(k+2)^2a_1^2,
\end{split}
\end{equation*}
Plugging into \eqref{Hamilton}, we have
\begin{equation*}
\begin{split}
a_1 &\geq a_1^2+b_1^2+2(a_2a_3+b_2b_3)\\
&> a_1^2+2ka_1[1-(k+1)a_1]-3(k+2)^2a_1^2\\
&=2ka_1 -(5k^2+14k+11)a_1^2,
\end{split}
\end{equation*}
therefore we get
\begin{equation}
\begin{split}
a_1 > \frac{2k-1}{5k^2+14k+11}.
\end{split} \label{solution1}
\end{equation}

On the other hand, by the Berger curvature decomposition, $|b_3-b_2|\leq a_3-a_2=1-(2k+1)a_1,$ so we have
\begin{equation} \label{b-term}
b_1^2+2b_2b_3=\frac{3}{2}b_1^2-\frac{1}{2}(b_3-b_2)^2\geq-\frac{1}{2}(a_3-a_2)^2\geq -\frac{1}{2}[1-(2k+1)a_1]^2,
\end{equation}
therefore,
\begin{equation*}
\begin{split}
a_1 &\geq a_1^2+b_1^2+2(a_2a_3+b_2b_3)\\
&\geq a_1^2+2ka_1[1-(k+1)a_1]-\frac{1}{2}[1-(2k+1)a_1]^2\\
&= -(4k^2+4k-\frac12)a_1^2+(4k+1)a_1-\frac{1}{2},
\end{split}
\end{equation*}
which implies
\begin{equation}
\begin{split}
a_1 &\leq \frac{4k-\sqrt{8k^2-8k+1}}{8k^2+8k-1} \hspace{0.5cm} \text{or}\hspace{0.5cm}  a_1 \geq \frac{4k+\sqrt{8k^2-8k+1}}{8k^2+8k-1}.
\end{split} \label{solution2}
\end{equation}

If $a_1 \geq \frac{4k+\sqrt{8k^2-8k+1}}{8k^2+8k-1}$, then $a_1=\frac{1}{3}$ if $k=1$; and if $k>1$ direct computation shows that,
\begin{equation*}
\begin{split}
a_2-a_3=(2k+1)a_1-1\geq (2k+1)\frac{4k+\sqrt{8k^2-8k+1}}{8k^2+8k-1}-1 &> 0,
\end{split}
\end{equation*}
which contradicts to $a_2\leq a_3$. Therefore from (\ref{solution1}) and (\ref{solution2}), we have either $a_1=\frac{1}{3}$, or
\begin{equation*}
\begin{split}
\frac{2k-1}{5k^2+14k+11} &< a_1 \leq \frac{4k-\sqrt{8k^2-8k+1}}{8k^2+8k-1},
\end{split}
\end{equation*}
which holds only if $1\leq k\leq 4$, so we get
\begin{equation*}
\begin{split}
a_1> \min_{1\leq k\leq 4}  \frac{2k-1}{5k^2+14k+11} =\frac{1}{30}.
\end{split}
\end{equation*}

\vspace{0.3cm}

\textit{At last we prove $a_1+a_2> 0$ implies $\mathfrak{R}$ is 4-positive.} It suffices to prove that $a_1+a_2> 0$ implies $1+(a_1\pm b_1)> 0$. From the Berger decomposition we have $|b_1|\leq \frac{1}{3}-a_1$, so $a_1>-\frac{1}{3}$ implies $1+(a_1\pm b_1)>0$. We will show that in fact $a_1+a_2>0$ implies $a_1>4-\sqrt{17}$.

Assuming $a_1(p)=\min a_1$. Plugging \eqref{b-term} into \eqref{Hamilton}, we have
\begin{equation} \label{eqn2}
\begin{split}
a_1(p)\geq& a_1^2+b_1^2+2(a_2a_3+b_2b_3)\\
\geq& a_1^2+2a_2a_3-\frac{1}{2}(a_3-a_2)^2\end{split}
\end{equation}
Since $a_3+a_2=1-a_1$, and $a_2>-a_1, a_3<1$, we have (the minimum is achieved on the boundary)
\begin{equation} \label{eqn3}
\begin{split}
2a_2a_3-\frac{1}{2}(a_3-a_2)^2=-\frac{1}{2}a_2^2-\frac{1}{2}a_3^2+3a_2a_3
>-\frac{1}{2}a_1^2-\frac{1}{2}-3a_1,
\end{split}
\end{equation}
Plugging \eqref{eqn3} into \eqref{eqn2}, we get that $a_1>4-\sqrt{17}$.
\qed

\begin{remark} In the author's thesis \cite{Wuthesis}, there was a naive mistake that ``by Berger curvature decomposition $a_1+a_2>0$ automatically implies $1+(a_1\pm b_1)>0$". The author caught and corrected this (see the last step in the proof of Theorem \ref{thm}) in August 2012 when he arrived at Cornell University as a postdoctoral fellow and prepared for seminar talks on his thesis and the work of Gursky and LeBrun \cite{GL} and Yang \cite{Yang}.
\end{remark}

The proof of Proposition \ref{remark1} contains a two-step constrained optimization. We omit the details since the main argument is the same as the proof of Theorem \ref{thm}.

Step one, we show that $K<\frac{14-\sqrt{19}}{12}$ implies $K>\frac{5-\sqrt{19}}{12}\approx0.0534$. Recall that at the minimum point of the sectional curvature, $a_1^2+b_1^2+2(a_2a_3+b_2b_3)\leq a_1$. Therefore the constrained optimization is
\begin{equation*}
\begin{split}
\textit{Minimize}&\qquad a_1,\\
\textit{Subject\ to}&\qquad a_3<\frac{14-\sqrt{19}}{12},\\
&\qquad a_1^2+b_1^2+2(a_2a_3+b_2b_3)\leq a_1,\\
&\qquad a_1+b_1 \leq a_2+b_2 \leq a_3+b_3,\\
&\qquad a_1-b_1 \leq a_2-b_2 \leq a_3-b_3,\\
&\qquad a_1+a_2+a_3=1,\ b_1+b_2+b_3=0.
\end{split}
\end{equation*}

Step two, we show that $K<\frac{14-\sqrt{19}}{12}$ \textit{and} $K>\frac{5-\sqrt{19}}{12}$ imply 3-positive curvature operator. To do this, we evaluate Equation (\ref{Hamiltonequation}) at eigenvectors and plug in the Berger decomposition. We denote eigenvalues of $\mathfrak{R}^+$ and $\mathfrak{R}^-$ by $\lambda_i=a_i+b_i$, $\mu_i=a_i-b_i$, and denote corresponding eigenvectors by $\omega_i^+$, $\omega_i^-$, respectively. We get
\begin{equation*}
\begin{cases}
\Delta R(\omega_1^+,\omega_1^+)+\lambda_1^2+2\lambda_2\lambda_3=\lambda_1,&\\
\Delta R(\omega_2^+,\omega_2^+)+\lambda_2^2+2\lambda_1\lambda_3=\lambda_2,&\\
\Delta R(\omega_3^+,\omega_3^+)+\lambda_3^2+2\lambda_1\lambda_2=\lambda_3,&\\
\Delta R(\omega_1^-,\omega_1^-)+\mu_1^2+2\mu_2\mu_3=\mu_1.\\
\Delta R(\omega_2^-,\omega_2^-)+\mu_2^2+2\mu_1\mu_3=\mu_2.\\
\Delta R(\omega_3^-,\omega_3^-)+\mu_3^2+2\mu_1\mu_2=\mu_3.
\end{cases}
\end{equation*}

Suppose the minimum of sum of any three eigenvalues is achieved by $\lambda_1+\lambda_2+\mu_1=1-\lambda_3+\mu_1=\min(I-\mathfrak{R}^++\mathfrak{R}^-)$ at a point $q$. Then at $q$, taking the sum we get
\begin{equation}\label{minimumthreeeigenvalue}
\boxed{\mu_1^2+2\mu_2\mu_3-\lambda_3^2-2\lambda_1\lambda_2\leq\mu_1-\lambda_3.}
\end{equation}
Therefore the constrained optimization is
\begin{equation*}
\begin{split}
\textit{Minimize}&\qquad 1+\mu_1-\lambda_3,\\
\textit{Subject\ to}&\qquad \lambda_3+\mu_3<\frac{14-\sqrt{19}}{6},\\
&\qquad \lambda_1+\mu_1>\frac{5-\sqrt{19}}{6},\\
&\qquad \mu_1^2+2\mu_2\mu_3-\lambda_3^2-2\lambda_1\lambda_2\leq\mu_1-\lambda_3.\\
&\qquad \lambda_1\leq\lambda_2\leq\lambda_3,\ \mu_1\leq\mu_2\leq\mu_3,\\
&\qquad \lambda_1+\lambda_2+\lambda_3=1,\ \mu_1+\mu_2+\mu_3=1.
\end{split}
\end{equation*}
We get $(1-\lambda_3+\mu_1)(q)>0$. If the minimum is attained by $\lambda_1+\mu_1+\mu_2$ at some point, then we get the same conclusion.
\qed

The proof of Proposition \ref{prop2} follows from an observation that at the minimum point of $\lambda_1$, one has $\lambda_1^2+2\lambda_2\lambda_3\leq\lambda_1$.
\qed


\begin{thebibliography}{20}
\bibliographystyle{alpha}
\renewcommand{\baselinestretch}{0}
\footnotesize

\bibitem{berger61} Berger, M., \emph{Sur quelques vari\'et\'es d'Einstein compactes}, Ann. Mat. Pura Appl. \textbf{53} (1961), 89--95.

\bibitem{Besse} Besse, A., \emph{Einstein manifolds}, Berlin-Heidelberg, Springer-Verlag, 1987.

\bibitem{BW} B\"ohm, C., Wilking, B., \emph{Manifolds with positive curvature operators are space forms},
Ann. of Math. (2) \textbf{167} (2008), 1079--1097.

\bibitem{Brendle} Brendle, S., \emph{Einstein manifolds with nonnegative isotropic curvature are locally symmetric},
Duke Math. Journal \textbf{151} (2010), 1--21.

\bibitem{Chen} Chen, H., \emph{Pointwise $\frac{1}{4}$-pinched 4-manifolds}, Ann. Global Geom. \textbf{9} (1991), 161--176.

\bibitem{Costa} Costa, \'E, \emph {On Einstein four-Manifolds}, J. of Geometry and Physics \textbf{51} (2004), 244--255.

\bibitem{Der} Derdzi\'nski, A., \emph{Self-dual K\"ahler manifolds and Einstein manifolds of dimension four}, Compositio Math. \textbf{49} (1983), 405--433.

\bibitem{GL} Gursky, M., LeBrun, C., \emph{On Einstein manifolds of positive sectional curvature}, Ann. Glob. An. Geom. \textbf{17} (1999), 315--328.

\bibitem{H3} Hamilton, R., \emph{Three-manifolds with positive Ricci curvature}, J. Differential Geometry \textbf{17} (1982), 255--306.

\bibitem{H4} Hamilton, R., \emph{Four-manifolds with positive curvature operator}, J. Differential Geometry \textbf{24} (1986), 153--179.

\bibitem{MM} Micallef M., Moore J. D., \emph{Minimal two-spheres and the topology of manifolds with positive curvature on totally isotropic two-planes},
Ann. of Math. \textbf{127} (1988), 199--227.

\bibitem{ST} Singer, I.M., A. Thorpe, J.A., \emph{The curvature of 4-dimensional Einstein spaces}, In Global Analysis (Papers in Honor
of K. Kodaira), 355--365. Univ. Tokyo Press, Tokyo, 1969.

\bibitem{Tachi} Tachibana, S., \emph{A theorem of Riemannian manifolds of positive curvature operator,} Proc. Japan Acad. \textbf{50} (1974), 301--302.

\bibitem{Wuthesis} Wu, P., \emph{Studies on Einstein manifolds and gradient Ricci solitons}, Ph. D. dissertation,
University of California, Santa Barbara, 2012.

\bibitem{Wu1} Wu, P., \emph{Einstein four-manifolds of three-nonnegative curvature operator}, preprint, 2013.

\bibitem{Wu2} Wu, P., \emph{A Weitzenbock formula for canonical metrics on four-manifolds}, Trans. Amer. Math. Soc. \textbf{369} (2017), 1079--1096.

\bibitem{Yang} Yang, D., \emph{Rigidity of Einstein 4-manifolds with positive curvature}, Invent. Math. \textbf{142} (2000), 435--450.

\end{thebibliography}
\end{document}